\documentclass[review,3p]{elsarticle}

\usepackage{amsmath}
\usepackage{amssymb}
\usepackage{hyperref}
\hypersetup{colorlinks = true,linkcolor = blue,anchorcolor =red,citecolor = blue,filecolor = red,urlcolor = red,
            pdfauthor=author}
\journal{}

\newdefinition{definition}{Definition}
\newdefinition{remark}{Remark}
\newdefinition{example}{Example}

\newtheorem{theorem}{Theorem}
\newtheorem{lemma}[theorem]{Lemma}

\newproof{proofOfTheoremCylinderStructure}{Proof of Theorem \ref{theorem:cylinder}}
\newproof{proofOfTheoremNonNegativeScalarCurvature}{Proof of Theorem \ref{theorem:nonnegative-scalar-curvature}}
\newproof{proofOfTheoremConcentratedLength}{Proof of Theorem \ref{theorem:concentrated-length}}
\newproof{proofOfTheoremCompactRiemannSolitonsAreTrivial}{Proof of Theorem \ref{theorem:compact-riemann-solitons-are-trivial}}
\newproof{proofOfLemmaWarpedProductStructure}{Proof of Lemma \ref{lemma:warped-product-structure}}
\newproof{proofOfLemmaConsequencesOfTheFundamentalEquation}{Proof of Lemma \ref{lemma:consequences-of-the-fundamental-equation}}
\newproof{proofOfLemmaAlmostSchoutenSoliton}{Proof of Lemma \ref{lemma:almost-schouten-soliton}}

\begin{document}

\begin{frontmatter}
\title{Rigidity results for Riemann solitons}

\author[1]{Willian Tokura\corref{cor1}\fnref{orcid1}}
\cortext[cor1]{Corresponding author}
\ead{willian.tokura@ufac.br}

\author[2]{Marcelo Barboza\fnref{orcid2}}
\ead{bezerra@ufg.br}

\author[2]{Elismar Batista\fnref{orcid3}}
\ead{elismardb@gmail.com}

\author[2]{Ilton Menezes\fnref{orcid4}}
\ead{iltomenezesufg@gmail.com}

\address[1]{Universidade Federal do Acre, CCET, 69920-900, Rio Branco - AC, Brazil}
\address[2]{Universidade Federal de Goi\'{a}s, IME, 74001-970, Goi\^{a}nia - GO, Brazil}

\fntext[orcid1]{https://orcid.org/0000-0001-9363-793X}
\fntext[orcid2]{https://orcid.org/0000-0003-4258-4676}
\fntext[orcid3]{https://orcid.org/0000-0002-4391-4238}
\fntext[orcid4]{https://orcid.org/0000-0002-9590-6731}

\begin{abstract}
  In this paper we take a look at conditions that make a Riemann soliton
  trivial, compacity being one of them. We also show that the behaviour at
  infinity of the gradient field of a non-compact gradient Riemann soliton might
  cause the soliton to be an Einstein manifold. Finally, we obtain scalar
  curvature estimates for complete shrinking or steady gradient Riemann solitons
  whose scalar curvature is bounded from below at infinity.
\end{abstract}

\end{frontmatter}

\section{Introduction and main results}
\label{sec:introduction}

The Riemann flow, which was introduced by Udriste \cite{udriste2010riemann,
udriste2011riemann}, is the flow associated with the evolution equation
\begin{equation}
  \label{eq:riemann-flow}
	\frac{\partial}{\partial t}G(t)=-2Rm_{g(t)},
\end{equation}
for Riemannian metrics on a given differentiable manifold $M^{n}$, where
$Rm_{g(t)}$ is the Riemann curvature $(0,4)$-tensor of the metric $g(t)$,
$G=\frac{1}{2}g\odot g$ and $\odot$ is the Kulkarni-Nomizu product. Recall
that, for symmetric $(0,2)$-tensors $A$ and $B$, the Kukarni-Nomizu product
$A\odot B$ is defined by (see \cite{besse2007einstein}):
\begin{align*}
	(A\odot B)(X_{1},X_{2},X_{3},X_{4}) & =
	A(X_{1},X_{3})B(X_{2},X_{4})+
	A(X_{2},X_{4})B(X_{1},X_{3})                                              \\
	                                    & \qquad-A(X_{1},X_{4})B(X_{2},X_{3})-
	A(X_{2},X_{3})B(X_{1},X_{4}).
\end{align*}

Short-time existence and uniqueness of the flow \eqref{eq:riemann-flow} is
guaranteed by Hirica and Udriste on compact manifolds in \cite{hirica2012basic},
where the authors utilize the Ricci flow to construct short-time solutions of
\eqref{eq:riemann-flow}.

Just like it happens with Ricci solitons, critical metrics for the Riemann flow
are self similar solutions of \eqref{eq:riemann-flow}, that is, it evolves overs
time from a given Riemannian metric on the manifold by diffeomorphisms and
dilatations \cite{hirica2016ricci}.

\begin{definition}
  A Riemannian manifold $(M^{n},g)$, $n\geqslant3$, is a Riemann soliton if
  there exist a vector field $X\in\mathfrak{X}(M)$ and a real number
  $\lambda\in\mathbb{R}$ satisfying:
	\begin{equation}
    \label{eq:riemman-soliton}
		Rm+\frac{1}{2}\mathcal{L}_{X}g\odot g=\lambda G,
	\end{equation}
  where $\mathcal{L}_{X}g$ is the Lie derivative of $g$ in the direction of $X$.
  We write the soliton in \eqref{eq:riemman-soliton} as $(M^{n},g,X,\lambda)$
  for the sake of simplicity. It may happen that $X=\nabla f$ for some $f\in
  C^{\infty}(M)$, in which case we say that $(M^n,g,\nabla f,\lambda)$ is a
  gradient Riemann soliton. Notice that equation \eqref{eq:riemman-soliton}
  for gradient Riemann solitons reads:
	\begin{equation}
    \label{eq:gradient-riemann-soliton}
		Rm+\nabla^2{f}\odot g=\lambda G.
	\end{equation}
  A Riemann soliton is classified into one of three types according to the sign
  of $\lambda$, expanding if $\lambda < 0$, steady if $\lambda = 0$ and,
  finally, shrinking if $\lambda > 0$. Also, a Riemann soliton is called trivial
  if either \(\mathcal{L}_{X}g=0\) or \(\nabla^{2}f=0\).
\end{definition}

One might start picturing Riemann solitons by noticing that the trivial ones
satisfy $Rm=\frac{\lambda}{2}g\odot g$, thus having constant sectional
curvature. Also, in \cite{catino2019potential} Catino and Mastrolia define what
they call \(X\)-space and \(f\)-space forms. Interestingly, \(X\)-space forms
with constant \(div(X)\) are Riemann solitons.

\begin{example}
  Let $(\mathbb{Q}^{n}(c),g)$ be a Riemannian manifold of constant sectional
  curvature $c\in\{-1,0,1\}$. Then, $(\mathbb{Q}^{n}(c),g,f,c)$ is a trivial
  Riemann soliton for any constant function \(f\) on \(\mathbb{Q}^{n}(c)\).
\end{example}

\begin{example}
  \textup{(Gaussian soliton)}
  The Gaussian soliton on $\mathbb{R}^{n}$, which is given by
  \[
    g_{ij}=\delta_{ij}
    \quad\mbox{and}\quad
    f(x)=\frac{\lambda}{4}|x|^{2},
  \]
  is a Riemann soliton since
  \[
    Rm=0
    \quad\mbox{and}\quad
    \nabla^{2}f=\frac{\lambda}{2}g.
  \]
\end{example}

\begin{example}\label{ex1}
  Let \((M^{n+1},g)\) be the Riemannian product manifold of \((0,\infty)\) and
  \(\mathbb{S}^{n}\). Then, \((M^{n+1},g,\nabla f,\lambda)\) is a gradient
  Riemann soliton for
  \[
    f:(0,\infty)\times\mathbb{S}^{n}\longrightarrow\mathbb{R},
    \quad
    (r,p)\longmapsto\frac{\lambda}{4r},
  \]
  where \(\lambda\in\mathbb{R}\).
\end{example}

\begin{example}\label{ex2}
  Let \((M^{n+1},g)\) be the product manifold \((0,\infty)\times\mathbb{S}^{n}\)
  equipped with the metric tensor given by
  \[
    g_{(r,p)}(t_{1}\oplus v_{1},t_{2}\oplus v_{2})=
    r
    \left(
      t_{1}t_{2}+v_{1}v_{2}
    \right),
  \]
  for every \((r,p)\in(0,\infty)\times\mathbb{S}^{n}\) and
  \(
    t_{1}\oplus v_{1},t_{2}\oplus v_{2}
    \in
    \mathbb{R}\oplus T_{p}\mathbb{S}^{n}
    \cong
    T_{(r,p)}\left((0,\infty)\times\mathbb{S}^{n}\right)
  \). Then, the quadruple \((M^{n+1},g,\nabla f,\lambda)\) is a gradient Riemann
  soliton for
  \[
    f:(0,\infty)\times\mathbb{S}^{n}\longrightarrow\mathbb{R},
    \quad
    (r,p)\longmapsto\ln{\left(r^{2}\right)}+\frac{r^{2}}{4}.
  \]
\end{example}

Examples \eqref{ex1} and \eqref{ex2} both are locally conformally flat gradient Riemann
solitons invariant by the left action of
\[
  \mbox{SO}(n,\mathbb{R})=
  \left\{
    A\in\mbox{GL}(n,\mathbb{R}):A^{-1}=A^{t}
  \right\},
\]
on \((0,\infty)\times\mathbb{S}^{n}\), given by
\[
  \mbox{SO}(n,\mathbb{R})\times\left((0,\infty)\times\mathbb{S}^{n}\right)
  \longrightarrow(0,\infty)\times\mathbb{S}^{n},
  \quad
  (A,(r,p))\longmapsto(r,A\cdot p).
\]
We refer the reader to \cite{barbosa2014gradient,batista2019warped,leandro2020reduction,neto2018gradient,pina2020gradient,tokura2019warped,tokura2021invariant}
for examples of other symmetric solitons. When we first met Riemann solitons it
occurred to us the results of Ivey \cite{ivey1993ricci}, showing that expanding
or steady compact gradient Ricci solitons are rigid, of Hsu \cite{hsu2012note},
showing that compact gradient Yamabe solitons are trivial, of Barros \textit{et
al.} \cite{barros2014compact}, showing that compact almost Ricci solitons with
constant scalar curvature are Einstein manifolds, of Catino \textit{et al.}
\cite{catino2016gradient}, showing that compact $\rho$-Einstein solitons are
trivial, and even one of our own, Tokura \textit{et al.}
\cite{tokura2021triviality}, showing that compact gradient $k$-Yamabe solitons
are trivial. On Riemann solitons, specifically, we found a result due to Blaga
\cite{blaga2021remarks}, showing that compact Riemann solitons with potential
vector field $X$ of constant length are trivial. Our first result extends
Blaga's result as it doesn't impose any additional conditions on the potential
vector field.

\begin{theorem}\label{theorem:compact-riemann-solitons-are-trivial}
	Compact Riemann solitons  $(M^{n},g,X,\lambda)$ are trivial.
\end{theorem}

Also, we show that non-compact Riemann solitons are Einstein manifolds as long
as it has constant scalar curvature and a certain weighted integral of
\(|\nabla f|^{2}\) is finite. 

\begin{theorem}
  \label{theorem:concentrated-length}
  Let $(M^{n},g,\nabla f,\lambda)$ be a complete non-compact gradient Riemann
  soliton with constant scalar curvature satisfying
	\begin{equation*}
		\int_{M^{n}\setminus B_{r}(x_{0})}
      \frac{|\nabla f|^{2}}{d(x,x_{0})^{2}}
		\,dv_{g}
    <\infty,
	\end{equation*}
  for some \(r>0\), where $d$ is the \(g\)-distance function on \(M^{n}\) and
  $B_{r}(x_{0})$ then is the \(d\)-ball of radius $r$ centered at $x_{0}$.
  Then, $(M^{n},g)$ is an Einstein manifold.
\end{theorem}

Next, we show that for Riemann solitons the sign os the scalar curvature at
infinity determines the sign of the scalar curvature as a whole.

\begin{theorem}
  \label{theorem:nonnegative-scalar-curvature}
  Let $(M^{n},g,X,\lambda)$ be a complete non-compact shrinking or steady
  Riemann soliton. Assume that $\liminf_{x\to\infty}S(x)\geqslant0$.
  Then, $M$ has nonnegative scalar curvature. Furthermore, if $M$ is not scalar
  flat, then $S>n\lambda$.
\end{theorem}

Our last result classifies scalar flat and steady gradient Riemann solitons.

\begin{theorem}\label{theorem:cylinder}
  Let $(M^{n},g,\nabla f,\lambda)$ be a complete steady gradient
  Riemann soliton. If $M^{n}$ is scalar flat and $f$ is not constant, then $M^{n}$ is isometric to a
  cylinder $\mathbb{R}\times \Sigma^{n-1} $ where $\Sigma^{n-1}$ is a Ricci flat
  hypersurface of $M^{n}$ and $f(t,x)=at+b$, for certain $a, b\in\mathbb{R}$,
  $a\neq0$.
\end{theorem}

\section{Key lemmas and proofs}

Our first lemma provides a link between Riemann solitons and almost Ricci
solitons. This relation was also noted by Hirica and Udriste
\cite{hirica2016ricci} and Blaga \cite{blaga2021remarks}.

\begin{lemma}\label{lemma:almost-schouten-soliton}
  The Riemannian manifold $(M^{n},g)$ admits a Riemannian soliton structure
  $(M^{n},g,X,\lambda)$ if and only if $(M^n,g)$ has null Weyl tensor and admits
  an almost Ricci soliton structure with soliton vector field $(n-2)X$
  satisfying 
	\begin{equation}\label{redução}
		Rc_{jk}+\frac{n-2}{2}(\mathcal{L}_{X}g)_{jk}=
		\left[
			\frac{(n-2)\lambda}{2}+\frac{S}{2(n-1)}
			\right]
		g_{jk},
	\end{equation}
  where \(Rc\) is the Ricci tensor and \(S\) the scalar curvature of
  $(M^{n},g)$.
\end{lemma}

\begin{proofOfLemmaAlmostSchoutenSoliton}
	From the fundamental equation \eqref{eq:riemman-soliton}, we have
	\begin{equation*}
		\begin{split}
			Rm_{ijkl}+
			\frac{1}{2}[g_{il}(\mathcal{L}_{X}g)_{jk}+
				g_{jk}(\mathcal{L}_{X}g)_{il}-
				g_{ik}(\mathcal{L}_{X}g)_{jl}-
				g_{jl}(\mathcal{L}_{X}g)_{ik}]=
			\lambda[g_{il}g_{jk}-g_{ik}g_{jl}],
		\end{split}
	\end{equation*}
	which by contraction over $i$ and $l$, gives
	\begin{equation}\label{equation:2}
		Rc_{jk}+\frac{n-2}{2}(\mathcal{L}_{X}g)_{jk}=
		[(n-1)\lambda-div(X)]g_{jk},
	\end{equation}
	or equivalently
	\begin{equation}\label{11}
		Rc_{jk}+\frac{n-2}{2}(\mathcal{L}_{X}g)_{jk}=
		\left[
			\frac{(n-2)\lambda}{2}+\frac{S}{2(n-1)}
			\right]
		g_{jk}.
	\end{equation}
	Now, substituting equation \eqref{11} into the decomposition of the Riemann
	tensor (see \cite{besse2007einstein}):
	\begin{equation}\label{Riemann1}
		Rm_{ijkl}=
		\left[
		\frac{1}{n-2}
		\left(
			Rc_{g}-\frac{S}{2(n-1)}g
			\right)
		\odot g
		\right]_{ijkl}+
		W_{ijkl},
	\end{equation}
	and using the Riemann soliton equation \eqref{eq:riemman-soliton}, we deduce
	that $W_{ijkl}=0$.
	
	On the other hand, if $(M^{n},g)$ has null Weyl $W=0$ and admits an almost Ricci soliton structure
	\[Rc_{jk}+\frac{n-2}{2}(\mathcal{L}_{X}g)_{jk}=
		\left[
			\frac{(n-2)\lambda}{2}+\frac{S}{2(n-1)}
			\right]
		g_{jk}.\]
Then from \eqref{Riemann1} we have
\begin{equation*}
		\begin{split}
		Rm_{ijkl}&=\left[
		\frac{1}{n-2}
		\left(
			Rc_{g}-\frac{S}{2(n-1)}g
			\right)
		\odot g
		\right]_{ijkl}=\frac{1}{2}\left[
		\left(
			-\mathcal{L}_{X}g+\lambda g
			\right)
		\odot g
		\right]_{ijkl}.
		\end{split}
	\end{equation*}
Therefore, $(M^n,g)$ admits a Riemann solitons structure
\[Rm_{ijkl}+\frac{1}{2}\left(\mathcal{L}_{X}g\odot g\right)_{ijkl}=\left(\frac{\lambda}{2}g\odot g\right)_{ijkl},\]
\end{proofOfLemmaAlmostSchoutenSoliton}

\begin{remark}In the particular case in which $X$ is a gradient vector field, we have that equation \eqref{redução} produces a Schouten soliton \cite{catino2016gradient}. Therefore, $(M^n,g)$ is a gradient Riemann soliton with potential $X$ if and only if $(M^n, g)$ is a Schouten soliton with potential $(n-2)X$ and $W_{ijkl}=0$.
\end{remark}
\begin{lemma}\label{lemma:consequences-of-the-fundamental-equation}
	Let $(M^{n},g,X,\lambda)$ be a Riemann soliton. Then the following formulas
	hold
	\begin{equation}\label{equation:3}
		div(X)=\frac{n\lambda}{2}-\frac{S}{2(n-1)}
	\end{equation}
	\begin{equation}\label{equation:4}
		Rc_{lj}X_{l}+\frac{S_{j}}{2(n-1)}+X_{jii}=0
	\end{equation}
	\begin{equation}\label{equation:5}
		\frac{n-2}{2} X_{j}S_{j}+
		\frac{n-2}{2}\lambda S+
		\frac{S^2}{2(n-1)}-|Rc|^2=0
	\end{equation}
\end{lemma}

\begin{proofOfLemmaConsequencesOfTheFundamentalEquation}
	In order to obtain \eqref{equation:3} it is enough to contract equation
	\eqref{redução}. Now, for equation \eqref{equation:4}, we remember the Ricci
	identity:
	\begin{equation}\label{www}
		X_{ijk}-X_{ikj}=X_{t}Rm_{tijk}.
	\end{equation}
	Contracting \eqref{www} with respect $i$, $k$, we have
	\[
		X_{iji}-X_{iij}=Rm_{liji}X_{l}=Rc_{lj}X_{l}.
	\]
	Next, differentiating  \eqref{redução} we get
	\begin{equation*}
		\begin{split}
			Rc_{ij,i}&=-\frac{n-2}{2}(X_{iji}+X_{jii})+\frac{S_{i}}{2(n-1)}g_{ij}\\
			&=-\frac{n-2}{2}(X_{iji}-X_{iij}+X_{iij}+X_{jii})+\frac{S_{i}}{2(n-1)}g_{ij}\\
			&=-\frac{n-2}{2}Rc_{lj}X_{l}-\frac{n-2}{2}(X_{iij}+X_{jii})+\frac{S_{i}}{2(n-1)}g_{ij}.
		\end{split}
	\end{equation*}
	Using the twice contracted second Bianchi identity and equation \eqref{equation:3}, we deduce
	\[
		\frac{1}{2}S_{j}=
		Rc_{ij,i}=
		-\frac{n-2}{2}Rc_{lj}X_{l}-
		\frac{n-2}{2}
		\left[
		-\frac{S_{j}}{2(n-1)}+X_{jii}
		\right]+
		\frac{S_{j}}{2(n-1)},
	\]
	which enables us obtain equation \eqref{equation:4}.

	Finally, from Lemma \ref{lemma:almost-schouten-soliton}, the Riemann soliton $(M^{n},g)$ admits a
	gradient almost Ricci solitons structure with vector field $(n-2)X$ and
	soliton function $\frac{(n-2)\lambda}{2}+\frac{S}{2(n-1)}$. Then, according to
	Lemma 3 of \cite{barros2014compact}, we have
	\begin{equation}\label{eqX1}
		\begin{split}
			\Delta_{(n-2)X}Rc_{ik}=
			\left[
				(n-2)\lambda+\frac{S}{(n-1)}
				\right]
			Rc_{ik}-2Rm_{ijks}Rc_{js}+
			\frac{n-2}{2}Rc_{is}(X_{sk}-
			X_{ks})\\
			+\frac{n-2}{2}Rc_{sk}(X_{si}-X_{is})+
			\frac{S_{ik}}{2}+\frac{\Delta S}{2(n-1)}g_{ki}-
			\frac{S_{ij}}{2(n-1)}g_{kj}.
		\end{split}
	\end{equation}
	Computing the trace of identity \eqref{eqX1}, we deduce
	\[
		-\frac{n-2}{2}X_{j}S_{j}=
		\left[
			\frac{n-2}{2}\lambda+\frac{S}{2(n-1)}
			\right]
		S-
		|Rc|^{2},
	\]
	which is equation \eqref{equation:5}.
\end{proofOfLemmaConsequencesOfTheFundamentalEquation}

\begin{proofOfTheoremCompactRiemannSolitonsAreTrivial}
	Let $p,q\in M$ be the points where $S$ attains its maximum and minimum in
	$M$, i.e.,
	\[
		S(q)=\min_{x\in M}S(x),\quad
		S(p)=\max_{x\in M}S(x),\quad
		\nabla S(p)=0=\nabla S(q),
	\]
	and
	\begin{equation*}
		\lambda-\frac{S(q)}{n(n-1)}\geqslant\lambda-
		\frac{S(x)}{n(n-1)}\geqslant\lambda-\frac{S(p)}{n(n-1)}.
	\end{equation*}
	From \eqref{equation:5} and the fact $|\textup{Rc}|^2\geqslant\frac{1}{n}S^2$,
	we have
	\begin{equation}\label{33}
		0=
		\frac{n-2}{2}\lambda S(q)+
		\frac{S(q)^2}{2(n-1)}-
		|\textup{Rc}|^2(q)\leqslant
		\frac{n-2}{2}S(q)
		\left(
		\lambda-
		\frac{S(q)}{n(n-1)}
		\right),
	\end{equation}
	and
	\begin{equation}\label{333}
		\begin{split}
			0=\frac{n-2}{2}\lambda S(p)+\frac{S(p)^2}{2(n-1)}-|\textup{Rc}|^2(p)
			\leqslant\frac{n-2}{2}S(p)\left(\lambda-\frac{S(p)}{n(n-1)}\right).
		\end{split}
	\end{equation}
	On the other hand, since $M$ is compact, integrating equation
	\eqref{equation:3}, we obtain the identity
	\begin{equation}\label{22}
		\lambda=\frac{1}{\textup{Vol}(M)}\int_{M}\frac{S}{n(n-1)}dv_{g},
	\end{equation}
	from what we see that
	\begin{equation}\label{Qe4}
		\lambda-\frac{S(q)}{n(n-1)}\geqslant0,\qquad
		\lambda-\frac{S(p)}{n(n-1)}\leqslant0.
	\end{equation}
	Now, we state that
	\[
		\lambda-\frac{S(q)}{n(n-1)}=0.
	\]
	In fact, if $\lambda-\frac{S(q)}{n(n-1)}>0$, then from \eqref{33}, we have
	$S(q)\geqslant0$ and $\lambda>0$. Consequently, $S(p)>0$ from
	\eqref{Qe4}. Now, since $S(p)>0$ we deduce from \eqref{333} and
	\eqref{Qe4} that $n(n-1)\lambda=S(p)$ and $S\equiv S(p)$ which is an
	absurd from $\lambda-\frac{S(q)}{n(n-1)}>0$. Therefore,
	$\lambda-\frac{S(q)}{n(n-1)}=0$.  However, this implies $S\equiv S(q)$ and
	$\lambda-\frac{S}{n(n-1)}=0$. Hence, from \eqref{redução} we deduce that
	$(M^{n},g)$ is a Ricci soliton.  From Perelman work's
	\cite{perelman2002entropy} we known that any compact Ricci soliton is a
	gradient Ricci soliton. Then there exist a smooth function $f$ on $M$ such
	that $div (X)=\Delta f=0$. The compacteness of $M$ implies that $f$ is
	constant and the soliton is trivial.
\end{proofOfTheoremCompactRiemannSolitonsAreTrivial}

\begin{proofOfTheoremConcentratedLength}
	We follow the argument in \cite{ma2012remarks}. The Bochner formula states
	that (\cite{bochner1954curvature}):
	\begin{equation*}
		\frac{1}{2}\Delta|\nabla f|^2=
		|\nabla^{2}f|^{2}+
		\langle
		\nabla f,\nabla\Delta f
		\rangle+
		Rc(\nabla f, \nabla f).
	\end{equation*}
	From equation \eqref{equation:4} we conclude that
	\begin{equation}\label{eeqq1}
		\begin{split}
			Rc_{lj}f_{l}&= -\frac{S_{j}}{2(n-1)}-(\nabla^{2} f)_{ji,i}\\
			&=-\frac{S_{j}}{2(n-1)}+\frac{Rc_{ji,i}}{n-2}-\frac{ S_{j}}{2(n-1)(n-2)}\\
			&=-\frac{S_{j}}{2(n-1)}+\frac{S_{j}}{2(n-2)}-\frac{S_{j}}{2(n-1)(n-2)}\\
			&=0.
		\end{split}
	\end{equation}
	Since $S$ is constant, we have that $div(X)=\Delta f$ is constant from
	\eqref{equation:3}. Hence, from \eqref{eeqq1}, we arrive at
	\begin{equation}\label{kkkk}
		\frac{1}{2}\Delta |\nabla f|^2=|\nabla^{2}f|^{2}.
	\end{equation}
	We choose  a cut-off function $\psi_{r}$ on $B_{2r}(x_{0})$ such that
	$0\leqslant \psi_{r}\leqslant 1$, $\mathrm{supp}(\psi_{r})\subset
		B_{2r}(x_{0})$ and
	\[
		\psi_{r}=1\quad\mbox{in}\quad B_{r}(x_{0}),\quad
		|\nabla\psi_{r}|^{2}\leqslant\frac{C}{r^{2}},\quad
		\Delta \psi_{r}\leqslant\frac{C}{r^{2}}.
	\]
	Then, multiplying \eqref{kkkk} by $\psi_{r}$ and integrating it over
	$B_{r}(x_{0})$, we get
	\begin{equation*}
		\begin{split}
			\int_{B_{r}(x_{0})}|\nabla^{2} f|^2\psi_{r}^{2}&=
			\int_{B_{r}(x_{0})}\frac{1}{2}\Delta |\nabla f|^2\psi_{r}^{2}\\
			&=\int_{B_{r}(x_{0})}\frac{1}{2}|\nabla f|^{2}\Delta\psi_{r}^{2}\\
			&\leqslant
			\int_{B_{2r}(x_{0})
				\setminus B_{r}(x_{0})}\frac{1}{2}\frac{C}{r^{2}}|\nabla f|^{2}
				\longrightarrow0\quad\mbox{as}\quad r\longrightarrow\infty.
		\end{split}
	\end{equation*}
	Hence
	\[
		\int_{M}|\nabla ^{2}f|^{2}=0.
	\]
	Therefore $\nabla^{2}f=0$, and from Lemma \ref{lemma:almost-schouten-soliton}
	$(M^{n},g)$ is Einstein.
\end{proofOfTheoremConcentratedLength}

\begin{proofOfTheoremNonNegativeScalarCurvature}
	From equation \eqref{equation:5} of Lemma
	\ref{lemma:consequences-of-the-fundamental-equation} and the inequality
	$|Rc|^{2}\geqslant\frac{S^{2}}{n}$, we deduce that
	\begin{equation}\label{label1}
		\langle X, \nabla S\rangle+\lambda S-\frac{S^2}{n(n-1)}\geqslant0.
	\end{equation}
	Using the maximum principle we can prove Theorem
	\ref{theorem:nonnegative-scalar-curvature}. In fact, suppose that
	$\inf_{M}S(x)<0$, then by the assumption
	$\liminf_{x\to\infty}S(x)\geqslant0$, there exists some point
	$z\in M$ such that
	\[
		S(z)=\inf_{M}S(x)<0.
	\]
	Hence,
	\[
		\Delta S(z)\geqslant0,\qquad \nabla S(z)=0.
	\]
	Then, from \eqref{label1} we deduce that
	\[
		\lambda S(z)-\frac{S(z)^{2}}{n(n-1)}\geqslant0.
	\]
	This is impossible since $\lambda S(z)-\frac{S(z)^{2}}{n(n-1)}<0$ for
	$\lambda\geqslant0$. Therefore $S(x)\geqslant0$. From the maximum principle we
	have that either $S\equiv0$ or $S(x)>n(n-1)\lambda$.
\end{proofOfTheoremNonNegativeScalarCurvature}

The next Lemma is a particular case of theorem 6.3 of Catino \textit{et al.}
\cite{catino2016geometry}. For the sake of completeness we will prove here.

\begin{lemma}\label{lemma:warped-product-structure}
	Let $(M^{n},g,\nabla f, \lambda)$ be a complete gradient Riemann soliton,
	$c\in\mathbb{R}$ be a regular value of $f$ and $\Sigma_{c}=f^{-1}(c)$ be its
	level surface. Then
	\begin{enumerate}[a)]
		\item
		      $|\nabla f|^{2}$ is constant on a connected component of $\Sigma_{c}$.
		      \vspace{0,1cm}
		\item
		      $\Sigma_{c}$ is totally umbilical and the scalar curvature $S$ is constant
		      on $\Sigma_{c}$
		      \vspace{0,1cm}
		\item
		      the mean curvature $H$ is constant on $\Sigma_{c}$.
		      \vspace{0,1cm}
		\item
		      $\Sigma_{c}$ is Einstein with respect to the induced metric.
		      \vspace{0,1cm}
		\item
		      in any open set of $\Sigma_{c}$ in which $f$ has non critical points, the
		      metric $g$ can be written as
		      \[
			      g=dr^{2}+\psi(r)^{2}g_{\Sigma},
		      \]
		      where $g_{\Sigma}$ is the metric induced by $g$ in $\Sigma_{c}$ and
		      $\psi(r)=e^{\frac{1}{n-1}\int_{r_{0}}^{r}H(s)ds}$.
	\end{enumerate}
\end{lemma}

\begin{proofOfLemmaWarpedProductStructure}
	Let $\{e_{1}, e_{2},\dots, e_{n}\}$ be a local orthonormal frame in a
	neighborhood of a regular point $p\in\Sigma_{c}$ such that $\{e_{2},\dots,
		e_{n}\}$ are tangent to $\Sigma_{c}$ and $e_{1}=\frac{\nabla f}{|\nabla f|}$.
	In order to prove $\textup{(a)}$, note that for any $X\bot\nabla f$ we have
	\[
		\nabla_{X}|\nabla f|^{2}=
		2\nabla^{2}f(X,\nabla f)=
		\frac{2}{n-2}Rc(X,\nabla f)=0.
	\]
	Now, from Lemma \ref{lemma:almost-schouten-soliton} the Weyl tensor $W$
	vanishes, then since Riemann solitons are in the class of almost Ricci
	solitons, we can apply Theorem 4.4 and Proposition 6.1 of
	\cite{catino2016geometry} with $\alpha=1$, $\beta=n-2$, $\mu=0$ to deduce that
	\begin{equation*}
		h_{ab}=\frac{H}{n-1}g_{ab},\qquad
		|\nabla S|=0,\quad\text{on}\quad\Sigma_{c},
	\end{equation*}
	where $h_{ab}$ and $H$ are respectively the second fundamental form and the
	mean curvature of $\Sigma_{c}$. This proves \textup{(b)}.

	For item \textup{(c)} we use the Codazzi equation
	\begin{equation}\label{Ti}
		Rm_{1cab}=
		\nabla^{\Sigma_{c}}_{a}h_{bc}-
		\nabla^{\Sigma_{c}}_{b}h_{ac},
		\qquad
		a, b, c \in\{2,\dots, n\}.
	\end{equation}
	Tracing over $b$ and $c$ in \eqref{Ti}, we obtain
	\begin{equation}\label{Ti2}
		Rc_{1a}=
		\nabla^{\Sigma_{c}}_{a}H-
		\nabla^{\Sigma_{c}}_{b}h_{ab}=
		\left(
		1-\frac{1}{n-1}
		\right)
		\nabla_{a}H.
	\end{equation}
	Then \textup{(c)} follows since $Rc_{1a}=0$.

	Now, using \eqref{redução} we deduce the following expression for the second
	fundamental form on $\Sigma_{c}$
	\begin{equation}\label{einstein2}
		h_{ab}=
		\langle
		\nabla_{a}e_{1},e_{b}
		\rangle=
		\frac{(\nabla^{2}f)_{ab}}{|\nabla f|}=
		\frac{
			\left[
				\frac{\lambda}{2}+\frac{S}{2(n-1)(n-2)}
				\right]
			g_{ab}-
			\frac{Rc_{ab}}{n-2}
		}{%
			|\nabla f|%
		}=
		\frac{H}{n-1}g_{ab}.
	\end{equation}
	Hence

	\begin{equation}\label{einstein3}
		Rc_{ab}=
		\left[
			\frac{\lambda(n-2)}{2}+\frac{S}{2(n-1)}-
			\frac{H(n-2)}{n-1}|\nabla f|
			\right]
		g_{ab}.
	\end{equation}
	Since $S$, $H$ and $|\nabla f|$ are constant on $\Sigma_{c}$, we have that
	$\Sigma_{c}$ is Einstein.

	Finally, since $\nabla f$ and the level surface of $f$ are orthogonal to each
	other, we can express the metric $g$ in the form
	\[
		g=dr^{2}+g_{ab}(r,\theta)d\theta^{a}d\theta^{b},
	\]
	where $\theta=(\theta^{2},\dots, \theta^{n})$ denotes a local coordinates on
	$\Sigma_{c_{0}}$ and $r(x)=\int\frac{df}{|\nabla f|}$. A good survey on level
	set structures can be found in \cite{cao2012locally, cao2012structure,
		catino2016geometry, leandro2021vanishing}.

	Fixing local coordinates system
	\[
		(x^{1},\dots, x^{n})=(r,\theta^{1},\cdots, \theta^{n}),
	\]
	we can express the second fundamental form of $\Sigma_{c}$ in terms of the
	Christoffel symbol $\Gamma_{ab}^{1}$, that is
	\[
		h_{ab}=
		-\langle
		\partial_{r},\nabla_{a}\partial_{b}
		\rangle=
		-\langle \partial_{r},\Gamma_{ab}^{1}\partial_{r}\rangle=
		-\Gamma_{ab}^{1}=
		-\frac{g^{11}}{2}
		\left(
		-\frac{\partial g_{ab}}{\partial r}
		\right)=
		\frac{1}{2}\frac{\partial g_{ab}}{\partial r}.
	\]
	Then
	\[
		\frac{1}{2}\frac{\partial g_{ab}}{\partial r}=\frac{H}{n-1}g_{ab},
	\]
	which implies that
	\[
		g_{ab}(r,\theta)=e^{\frac{2}{n-1}\int_{r_{0}}^{r}H(s)ds}g_{ab}(r_{0},\theta).
	\]
	Here the level surface $\{r = r_{0}\}$ corresponds to $\Sigma_{c_{0}}$.
	Therefore,
	\[
		g=
		dr^{2}+
		e^{\frac{2}{n-1}\int_{r_{0}}^{r}H(s)ds}
		g_{ab}(r_{0},\theta)d\theta^{a}d\theta^{b},
	\]
	which proves item \textup{(e)}.
\end{proofOfLemmaWarpedProductStructure}

\begin{proofOfTheoremCylinderStructure}
	From equation \eqref{equation:5} of Lemma
	\ref{lemma:consequences-of-the-fundamental-equation}, we deduce that $Rc=0$, so,  from
	\eqref{redução} we deduce that $\nabla^{2}f=0$. In particular $|\nabla f|$ is
	constant. Then either $f$ is constant, or $f$ has no critical point at all.
	Since $f$ is not constant, this latter case occurs. Now, Lemma
	\ref{lemma:warped-product-structure} implies that $(M^{n},g)$ is isometric to
	the warped product $\mathbb{R}\times_{\psi}\Sigma^{n-1}$ of the entire real
	line with a $(n-1)$–dimensional complete Riemannian manifold.

	We now prove that $\psi$ is constant and $\Sigma^{n-1}$ is Ricci flat. Indeed,
	since $\nabla^{2}f=0$, we have from equation \eqref{einstein2} that
	\begin{equation*}
		h_{ab}=
		\langle
		\nabla_{a}e_{1},e_{b}
		\rangle=
		\frac{(\nabla^{2}f)_{ab}}{|\nabla f|}=0,
	\end{equation*}
	which implies that $H=0$. Then
	\[
		\psi(r)^{2}=e^{\frac{2}{n-1}\int_{r_{0}}^{r}H(s)ds}=1.
	\]
	This fact shows that $M^{n}=\mathbb{R}\times_{1}\Sigma^{n-1}$ is a Riemannian
	product. On the other hand, from the expression \eqref{einstein3}, we have
	that $Rc_{ab}=0$, i.e., $\Sigma^{n-1}$ is Ricci flat. Furthermore, since
	$\nabla^{2}f=0$ on $\mathbb{R}\times \Sigma^{n-1}$, we derive that
	$f(t,x)=at+b$, $a,b\in\mathbb{R}$, $a\neq0$.
\end{proofOfTheoremCylinderStructure}


\section*{Conflict of interest}
The authors declare that there is no conflict of interest.

\section*{Data Availability Statements}
 Data sharing not applicable to this article as no datasets were generated or analysed during the current study.

\bibliography{main}

\end{document}